\documentclass[runningheads]{llncs}

\usepackage[T1]{fontenc}
\usepackage{graphicx}
\usepackage{mathtools}

\title{Reconfiguration Using Generalized Token Jumping}

\author{
  Jan Maty{\'a}{\v{s}} K{\v{r}}i{\v{s}}{\v{t}}an\inst{1}\orcidID{0000-0001-6657-0020} \and
  Jakub Svoboda\inst{2}
}

\authorrunning{J. M. K{\v r}i{\v s}{\v t}an and J. Svoboda}

\institute{
  Faculty of Information Technology, Czech Technical University in Prague, Czech Republic\\
  \email{kristja6@fit.cvut.cz} \and
  Institute of Science and Technology, Austria
}

\newtheorem{observation}{Observation}

\usepackage[utf8]{inputenc} 
\usepackage[nocompress]{cite}
\usepackage{todonotes}

\usepackage{amssymb} 
\usepackage{mathtools}
\usepackage{algorithm}
\usepackage{xspace}
\usepackage[noend]{algpseudocode}
\usepackage{thmtools} 
\usepackage{hyperref}
\usepackage{url}
\usepackage{float}

\usepackage{amsmath} 
\usepackage{thm-restate}
\usepackage{array}
\usepackage{tablefootnote}

\usepackage{tikz}
\usepackage{pgfplots}
\pgfplotsset{compat=newest}

\usepackage{cleveref}
\usepackage{subcaption}

\usetikzlibrary{positioning,arrows.meta}

\graphicspath{{./images/}}

\makeatletter
\newcommand{\customlabel}[2]{%
\protected@write \@auxout {}{\string \newlabel {#1}{{#2}{\thepage}{#2}{#1}{}} }%
\hypertarget{#1}{#2}
}
\makeatother

\newcommand{\decprob}[3]{
  \vspace{2mm}
\noindent\fbox{
  \begin{minipage}{0.96\textwidth}
  \textsc{#1}\\
    \textbf{Input:} #2 \\
  \textbf{Output:} #3
  \end{minipage}
  }
  \vspace{2mm}
}

\newcommand{\calO}{\mathcal{O}}
\newcommand{\TJ}{\textsc{Token Jumping}\xspace}
\newcommand{\TS}{\textsc{Token Sliding}\xspace}

\newcommand{\infTJ}{\textsc{$(k, 1)$-Token Jumping}\xspace}
\newcommand{\probVC}{\textsc{Vertex Cover}\xspace}
\newcommand{\probIS}{\textsc{Independent Set}\xspace}
\newcommand{\probDS}{\textsc{Dominating Set}\xspace}

\newcommand{\YES}{\textsf{YES}\xspace}
\newcommand{\NP}{\textsf{NP}\xspace}
\newcommand{\PSPACE}{\textsf{PSPACE}\xspace}
\newcommand{\NPSPACE}{\textsf{NPSPACE}\xspace}
\newcommand{\Pclass}{\textsf{P}\xspace}

\newcommand{\nkTJ}[2]{$(#1, #2)$-\textsc{Token Jumping}}
\newcommand{\nknkTJ}[4]{$\{(#1, #2), (#3, #4)\}$-\textsc{Token Jumping}}

\crefname{observation}{observation}{observations}  

\begin{document}

\maketitle

\begin{abstract}
  In reconfiguration, we are given two solutions to a graph problem, such as \probVC or \probDS, with each solution represented by a placement of tokens on vertices of the graph.
  Our task is to reconfigure one into the other using small steps while ensuring the intermediate configurations of tokens are also valid solutions.
  The two commonly studied settings are \TJ and \TS, which allows moving a single token to an arbitrary or an adjacent vertex, respectively.

  We introduce new rules that generalize \TJ, parameterized by the number of tokens allowed to move at once and by the maximum distance of each move.
  Our main contribution is identifying minimal rules that allow reconfiguring any possible given solution into any other for \probIS, \probVC, and \probDS.
  For each minimal rule, we also provide an efficient algorithm that finds a corresponding reconfiguration sequence.

  We further focus on the rule that allows each token to move to an adjacent vertex in a single step.
  This natural variant turns out to be the minimal rule that guarantees reconfigurability for \probVC.
  We determine the computational complexity of deciding whether a (shortest) reconfiguration sequence exists under this rule for the three studied problems.
  While reachability for \probVC is shown to be in \Pclass, finding a shortest sequence is shown to be \NP-complete.
  For \probIS and \probDS, even reachability is shown to be \PSPACE-complete.
\end{abstract}

\section{Introduction}

Reconfiguration problems arise whenever we want to transform one feasible solution into another in small steps while keeping all intermediate solutions also feasible.
The problem is widely studied in the context of graph problems, such as \textsc{Independent Set}~\cite{Lokshtanov_2019, Demaine_2015, Bonamy_2017, Belmonte_2020, Bartier_2021, Bonsma_2014, Yamada_2021, Hoang_2019}, \probVC~\cite{Mouawad_2018, ITO_2016}, \textsc{Dominating Set}~\cite{Bonamy_2021, Bousquet_2021, Lokshtanov_2018, Haddadan_2016, K_i_an_2023}, \textsc{Shortest Paths}~\cite{Kami_ski_2011, Gajjar_2022}, and \textsc{Coloring}~\cite{Bonsma_2009, Bonamy_2013, Bonsma_2014_2, Cereceda_2010}.
Moreover, outside of graphs, many results concern \textsc{Satisfiability}~\cite{Gopalan_2006, Mouawad_2017}.
See~\cite{Nishimura_2018} for a general survey.

In graph settings, the problem is characterized by a condition, a reconfiguration rule, an input graph, and an initial and a target configuration.
Configurations (including initial and target) are characterized by tokens placed on the vertices of the graph.
The condition restricts the placement of the tokens, for instance, the tokens need to form an independent set or a dominating set.
The reconfiguration rule allows only certain movements of the tokens to obtain a new configuration.
The goal is to move the tokens to the target configuration following the prescribed reconfiguration rule and ensuring all intermediate configurations satisfy the condition.

There are two main questions about reconfiguration problems.
Is the reconfiguration from the initial to the target configuration possible?
And if so, what is the smallest number of steps that reconfigure the initial to the target configuration?

Two reconfiguration rules are most commonly studied, \textsc{Token Jumping} and \textsc{Token Sliding}.
In \TJ, only one token is allowed to be removed from the graph and then placed on another vertex.
In the case of \TS, one token can be moved only to a neighboring vertex.
We can understand \TS as a restriction of \TJ to moves of distance $1$.


We generalize these reconfiguration rules.
First, we examine the case where all $k$ tokens can move by distance at most $d$ and we call this rule \nkTJ{k}{d}.
Second, we allow all $k$ tokens to move to distance $1$, and out of those, $k'$ tokens can move even farther to distance $d$.
We denote this rule as \nknkTJ{k}{1}{k'}{d}.
Note that \TS is equivalent to \nkTJ{1}{1} and on connected graphs the standard \textsc{Token Jumping} is equivalent to \nkTJ{1}{n}.

In this work, we present minimal rules that for a given condition (\probVC, \probDS, and \probIS) ensure that the reconfiguration problem is feasible for
any pair of initial and target configurations in a connected graph.
Conversely, we provide instances of initial and target configurations such that weaker rules do not allow reconfiguring one to the other.


We further consider the complexity of the reconfiguration problems under \nkTJ{k}{1} where $k$ is the size of the reconfigured solution, that is, each token can move to a neighbor.
We show that the reachability question, i.e.\ whether a reconfiguration sequence exists, is easy for \probVC under \infTJ.
On the other hand, deciding whether a reconfiguration sequence of length at most $\ell$ exists is \NP-hard, even for fixed $\ell \geq 2$.
Even just the reachability question for \probDS and \probIS is shown to be \PSPACE-hard.

It is known that reachability is \PSPACE-complete for \probIS under both \TJ and \TS even on graphs of bounded bandwidth~\cite{Wrochna_2018}.
The same result applies also for \probVC under \TJ and \TS if we do not allow multiple tokens on one vertex.
It is also known that reachability is \PSPACE-complete for \probDS under \TS even on bipartite graphs, split graphs, and bounded bandwidth graphs~\cite{Bonamy_2021}.
The same problem is also \PSPACE-complete under \TJ~\cite{Bonamy_2021}.

Reconfiguration by exchanging $k$ consecutive vertices at once has been investigated for reconfiguration of shortest paths~\cite{Gajjar_2022}.
Parallel reconfiguration has also been previously studied in the context of distributed computation~\cite{Censor_Hillel_2020}.
The rule of moving all tokens at once (or in other contexts, agents) appears in various other problems, such as in multi-agent path finding~\cite{Stern_2019} and graph protection problems~\cite{Klostermeyer_2016}.

In the context of multi-agent pathfinding, the \nkTJ{k}{1} rule is equivalent to the anonymous setting, meaning the agents are not distinguished.
A reconfiguration sequence under \nkTJ{k}{1} is thus equivalent to a solution to anonymous multi-agent pathfinding with an additional constraint on the positions of agents, such as requiring that the agents form a vertex cover or an independent set.
A certain case of such a global constraint has been studied, the condition where labeled agents must form a connected subgraph has received a lot of attention recently~\cite{Tateo_2018, Charrier_2020}.

\paragraph{Results.}
First, we explore which generalizations of \TJ allow us to guarantee that any two solutions can be reached from each other.
We show that if all $k$ tokens are not allowed to move at once, then we can not guarantee reachability for any of the three studied conditions (\probVC, \probDS, and \probIS).
That is, each of these problems requires at least \nkTJ{k}{1}.

Thus, we further focus on rules that allow all tokens to move at once and search for minimal possible extensions that guarantee reachability.
\Cref{tab:my-table} shows how much we have to extend the distance to which each token can move to ensure this guarantee.
We also provide corresponding examples where a weaker rule is not sufficient, thus the established results are tight with regard to the allowed distance.

\begin{table}[h]
    \centering
    \begin{tabular}{|l|p{3cm}|p{2.5cm}|p{2.5cm}|}
        \hline
        \textbf{Condition / Distance} & \textbf{$d = 1$} & \textbf{$d = 2$} & \textbf{$d \geq 3$} \\ \hline
        \probVC               & \multicolumn{3}{p{7.5cm}|}{Yes~(\Cref{thm:reconf-hg})} \\ \hline
        \probDS             & No (\Cref{prop:ds-not-connected}) & \multicolumn{2}{p{5cm}|}{Yes~(\Cref{prop:ds-connected})} \\ \hline
        \probIS            & \multicolumn{2}{p{5cm}|}{No~(\Cref{prop:indep-not-connected})} & Yes~(\Cref{prop:indep-possible}) \\ \hline
    \end{tabular}
    \caption{This table displays whether it is guaranteed that any solution can be reached from any other under \nkTJ{k}{d}, where $k$ is the size of the solution. That is, each token can move to distance at most $d$.}
    \label{tab:my-table}
\end{table}

We further explore if it suffices that only a part of the tokens moves by distance at most $d$, while the rest moves only by distance at most $1$.
For \probVC, the \nkTJ{k}{1} rule is already optimal, as any rule that allows moves of only $k - 1$ or fewer tokens has instances that can not be reconfigured, even if we allow any greater distance.

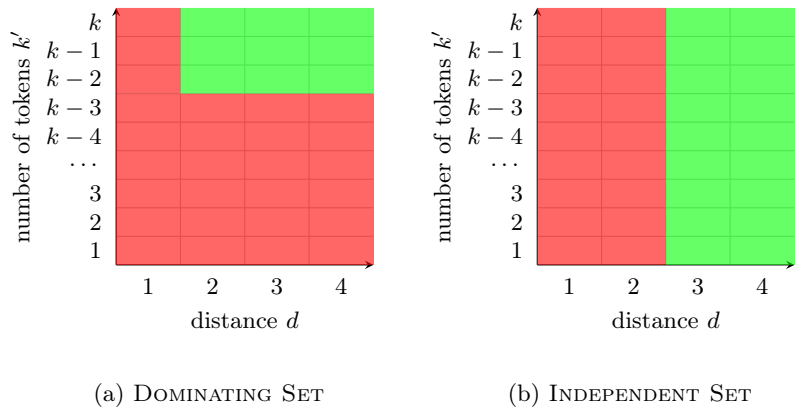
\begin{figure}
  \centering
  \begin{subfigure}{.45\textwidth}
    \centering
    \pgfplotsset{compat=newest}
\begin{tikzpicture}
    \begin{axis}[
    axis lines = left,
    width=5cm,
    height=5cm,
    ylabel = {number of tokens $k'$},
    xlabel = {distance $d$},
    xmin=1.5, xmax=5.5,
    ymin=1.5, ymax=10.5,
    xtick={0,1,2,3,4,5},
    ytick={0,1,2,3,4,5,6,7,8,9,10},
    xticklabel style={align=center},
    xticklabels={,0,1,2,3,4},
    yticklabel style={align=center},
    yticklabels={,0,1,2,3,$\dots$,$k-4$,$k-3$,$k-2$,$k-1$,$k$},
    minor tick num=1, 
    grid=minor, 
    tick style={draw=none}, 
    after end axis/.code={
        \path (axis cs:5,0) --
        (axis cs:6,0) node [anchor=west] {};
    }
    ]

    \addplot [
        draw=none,
        fill=red,
        fill opacity=0.6
    ] coordinates {
        (6.5,0) (0,0) (0,7.5) (6.5,7.5)
    } \closedcycle;

    \addplot [
        draw=none,
        fill=red,
        fill opacity=0.6
    ] coordinates {
        (2.5,7.5) (0,7.5) (0,10.5) (2.5,10.5)
    } \closedcycle;

    \addplot [
        draw=none,
        fill=green,
        fill opacity=0.6
    ] coordinates {
        (6.5,7.5) (2.5,7.5) (2.5,10.5) (6.5,10.5)
    } \closedcycle;

    \end{axis}
\end{tikzpicture}
    \captionsetup{justification=centering}
    \caption{\probDS}\label{fig:intervals}
  \end{subfigure}
  \begin{subfigure}{.45\textwidth}
    \centering
    \pgfplotsset{compat=newest}
\begin{tikzpicture}
    \begin{axis}[
    axis lines = left,
    width=5cm,
    height=5cm,
    ylabel = {number of tokens $k'$},
    xlabel = {distance $d$},
    xmin=1.5, xmax=5.5,
    ymin=1.5, ymax=10.5,
    xtick={0,1,2,3,4,5},
    ytick={0,1,2,3,4,5,6,7,8,9,10},
    xticklabel style={align=center},
    xticklabels={,0,1,2,3,4},
    yticklabel style={align=center},
    yticklabels={,0,1,2,3,$\dots$,$k-4$,$k-3$,$k-2$,$k-1$,$k$},
    minor tick num=1, 
    grid=minor, 
    tick style={draw=none}, 
    after end axis/.code={
        \path (axis cs:5,0) --
        (axis cs:6,0) node [anchor=west] {};
    }
    ]

    \addplot [
        draw=none,
        fill=red,
        fill opacity=0.6
    ] coordinates {
        (6.5,0) (0,0) (0,1.5) (6.5,1.5)
    } \closedcycle;

    \addplot [
        draw=none,
        fill=red,
        fill opacity=0.6
    ] coordinates {
        (3.5,1.5) (0,1.5) (0,10.5) (3.5,10.5)
    } \closedcycle;

    \addplot [
        draw=none,
        fill=green,
        fill opacity=0.6
    ] coordinates {
        (6.5,1.5) (3.5,1.5) (3.5,10.5) (6.5,10.5)
    } \closedcycle;

    \end{axis}
\end{tikzpicture}
    \captionsetup{justification=centering}
    \caption{\probIS}\label{fig:move2}
  \end{subfigure}
  \caption{Results for \nknkTJ{k}{1}{k'}{d}, where $k$ is the size of the reconfigured solution. Green shows where the reachability is guaranteed. Red means that reachability can not be guaranteed and a counterexample is shown.}\label{fig:is-ds-guarantees}
\end{figure}

The case of \probDS and \probIS is more complicated, the overview of results for \probDS and \probIS is shown in \Cref{fig:is-ds-guarantees}.
For \probDS, we establish that \nknkTJ{k}{1}{k - 3}{d} is not sufficient for all $d$.
On the other hand, we show that \nknkTJ{k}{1}{k-2}{2} is sufficient to guarantee reachability.
For \probIS, we show that \nknkTJ{k}{1}{1}{3} is sufficient and we show that \nkTJ{k}{2} is not sufficient for every solution size $k \geq 3$, i.e.\ allowing a move to distance $3$ is necessary.
For each of these reachability guarantees, we also provide a polynomial algorithm that finds a reconfiguration sequence.

When reachability is not guaranteed, it is natural to ask what is the complexity of the corresponding decision problem: are the two given solutions reachable under a given rule?
And if so, what is the minimum length of a reconfiguration sequence?
We establish the complexity of these two problems under \infTJ, see \Cref{tab:complexities} for an overview.

Note that deciding if any sequence exists can be reduced to deciding whether a sequence of length at most $2^n$ exists, as any reconfiguration sequence without repeating configurations will have at most $2^n$ moves.
Thus, if deciding whether a sequence exists is \PSPACE-complete, deciding if a sequence of length at most $\ell$ exists is \PSPACE-complete as well.

\begin{table}[h!]
  \centering
  \begin{tabular}{|r|l|l|}
    \hline
    \textbf{Condition / Deciding if exists} & Any sequence & Sequence with at most $\ell$ moves  \\
    \hline
    \probVC & \Pclass~(\Cref{thm:reconf-hg}) & \NP-complete~(\Cref{thm:vc-shortest-hard}) \\
    \hline
    \probIS & \multicolumn{2}{c|}{\PSPACE-complete~(\Cref{thm:pts-hard})} \\
    \hline
    \probDS & \multicolumn{2}{c|}{\PSPACE-complete~(\Cref{thm:ds-reconf})} \\
    \hline
  \end{tabular}
  \caption{
    Results on complexities of reconfiguration problems under \infTJ.
  }\label{tab:complexities}
\end{table}


\section{Preliminaries}



\paragraph{Graphs and graph problems.}
All graphs are assumed to be simple.
By $N(v)$ we denote the set of neighbors of $v$, and $N[v] = N(v) \cup \{v\}$ denotes the closed neighborhood of $v$.
By $dist_G(v, u)$ we denote the distance between $u$ and $v$ in $G$ and omit $G$ if it is clear from the context.
By $dist(C, u)$ where $C$ is a set of vertices, we denote the minimum distance from $u$ to any vertex of $C$.
Throughout the paper, we use $n$ to denote the number of vertices and $m$ to denote the number of edges of a graph or a hypergraph.
We use $\triangle$ to denote the symmetric difference of two sets, that is $A \triangle B = (A \setminus B) \cup (B \setminus A)$.

An \emph{independent set} is a set of vertices such that no two vertices are adjacent.
A \emph{vertex cover} is a set of vertices such that each edge contains at least one vertex in the set.
A \emph{dominating set} is a set of vertices $D$, such that each vertex is in $D$ or is incident to some vertex in $D$.

Let $H = (V, E)$ be a hypergraph, we define a path in $H$ as a sequence of vertices $v_1,v_2,\dots,v_\ell$ where $v_{i}$ and $v_{i+1}$ share an edge for all $1 \leq i \leq n - 1$ and the length of such path is $\ell - 1$.
The distance of two vertices in $H$ is the minimum length of a path between them.

\paragraph{Reconfiguration.}
Several previous works on \TS have allowed multiple tokens to occupy a single vertex~\cite{Bonamy_2021, K_i_an_2023, Bousquet_2021}.
We provide our results for the variant where tokens can not share vertices, which in case of upper bounds leads to stronger results.
At the same time, it is easy to verify that our lower bounds apply also for variants where multiple tokens can occupy the same vertex.

Let $\Pi(G)$ be the system of sets, which are feasible solutions of size $k$ of a given problem on $G$.
We say that sequence $D_1, \dots, D_\ell$ is a \emph{reconfiguration sequence} under the rule \nkTJ{k'}{d} if $D_i \in \Pi(G)$ for all $i$ and there is a \emph{move} under \nkTJ{k'}{d} from $D_{i}$ to $D_{i+1}$ for all $i < \ell$.

A move from $D_i$ to $D_{i+1}$ under \nkTJ{k'}{d} is a bijection $f : D_i \rightarrow D_{i+1}$ such that $f(v) \neq v$ for at most $k'$ distinct $v \in D_i$ and $dist_G(v, f(v)) \leq d$ for all $v$.
A move from $D_i$ to $D_{i+1}$ under \nknkTJ{k}{1}{k'}{d} is a bijection $f : D_i \rightarrow D_{i+1}$ such that there is $D' \subseteq D_i$ with $|D'| \leq k'$ such that $dist_G(v, f(v)) \leq 1$ for all $v \in D_i \setminus D'$ and $dist_G(v, f(v)) \leq d$ for all $v \in D'$.

Note that the standard definition of \TS coincides with \nkTJ{1}{1}, and the standard definition of \TJ coincides with \nkTJ{1}{n} on connected graphs.

We study the complexity of the two following problems.

\decprob{Reconfiguration of $\Pi$ under \nkTJ{k}{1}}
{Graph $G = (V, E)$ and two solutions $V_s, V_t \in \Pi(G)$ of size $k$.}
{Whether there exists a reconfiguration sequence between $V_s$ and $V_t$ with condition $\Pi$ under \nkTJ{k}{1}}

\decprob{Shortest reconfiguration of $\Pi$ under \nkTJ{k}{1}}
{Graph $G = (V, E)$, two solutions $V_s, V_t \in \Pi(G)$ of size $k$ and integer $\ell$.}
{Whether there exists a reconfiguration sequence between $V_s$ and $V_t$ of length at most $\ell$ with condition $\Pi$ under \nkTJ{k}{1}.}


\paragraph{Toolbox.}
Many of the results use the well-known Hall's Theorem, which we recall here.
A matching is a set of edges such that no two share an endpoint.
We say that a matching $M$ saturates $X$ if every vertex in $X$ is in some edge of $M$.
\begin{theorem}[Hall's Marriage Theorem \cite{Hall_1935}]\label{thm:hall}
Let \(G = (X \cup Y, E)\) be a bipartite graph.
A matching that saturates \(X\) exists if and only if for every subset \(S \subseteq X\), it holds \(|N(S)| \geq |S|\).
\end{theorem}

Recall that if a matching that saturates $X$ exists, it can be found in polynomial time by a well-known algorithm~\cite{west}.
If a matching saturating $X$ does not exist, it follows that there exists $S \subseteq X$ such that $|S| > |N(S)|$, such $S$ is called a \emph{Hall violator}.
A Hall violator (if exists) can be found in polynomial time by a well-known algorithm, well described for instance in~\cite{Gan_2019}.

Given a hypergraph $H$ and its two configurations of tokens $C_1, C_2$, finding a move or determining that no move exists under \nkTJ{k}{1} can be done efficiently by reducing to bipartite matching in the following auxiliary graph.
We define $B(C_1, C_2, H) = (V_1 \cup V_2, E)$, where the vertices $V_i = \{v^i_1, v^i_2, \dots, v^i_k\}$ correspond to vertices $C_i = \{v_1, v_2, \dots, v_k\}$ and
$\{v^1_x, v^2_y\} \in E$ if and only if $dist_H(v_x, v_y) \leq 1$.
\begin{observation}\label{obs:computing-move}
    Let $C_1$ and $C_2$ be two token configurations on a hypergraph $H$ and $d$ an integer.
    We can either find a move from $C_1$ to $C_2$ under $(k, 1)$-\TS or determine that no such move exists in $\mathcal{O}(n^{2.5} + n m)$ time.
\end{observation}
\begin{proof}
    There is a move from $C_1$ to $C_2$ if and only if $B(C_1, C_2, H)$ has a perfect matching.
    Indeed, if the perfect matching exists, it describes the move and every move can be described by a perfect matching.
    We can construct $B' = B(C_1, C_2, H)$ in time $\mathcal{O}(n^2 + nm)$ by iterating all edges of $H$.
    Note that $B'$ has $\mathcal{O}(n)$ vertices and $\mathcal{O}(n^2)$ edges its maximum matching can
    be found in $\mathcal{O}(n^{2.5})$~\cite{west} time.
\end{proof}


\section{Guaranteeing reconfigurability}

In this section, we examine under which reconfiguration rules are solutions to problems \probVC, \probDS, and \probIS always reconfigurable, and for which rules it may be impossible.
First, we show that $k-1$ tokens moving at once is not enough for any of these problems, even if we allow any greater distance.
Then, we describe a procedure that finds a reconfiguration sequence for vertex covers on hypergraphs and we use it to show that \nkTJ{k}{1} is sufficient for vertex cover.
With slight modification, a similar proof is used also for \probDS, where \nknkTJ{k}{1}{k-2}{2} rule is shown to be sufficient.
For that case, we also present a matching lower bound.
Finally, we describe a way how to reconfigure \probIS with \nknkTJ{k}{1}{1}{3} rule, and again, we show that this is tight.

\begin{proposition}\label{prop:need-k-to-move}
  For every integer $k \geq 2$ and problems \probVC, \probDS, \probIS, there exists a connected graph $G$ that has two unreachable solutions of size $k$ under \nkTJ{k - 1}{d} for every $d$.
\end{proposition}
\begin{proof}
  For \probVC and \probIS, consider a cycle of size $2k$.
  Such a cycle has exactly two minimum vertex covers of size $k$, both vertex disjoint.
  Thus, to reconfigure between them, each token has to move at once.
  Such a cycle also has exactly two maximum independent sets of size $k$, again both vertex disjoint.

  For \probDS, consider a cycle of size $3k$.
  Such a cycle has exactly three minimum dominating sets, all vertex disjoint.
\end{proof}

First, we make an observation regarding the case when no condition is placed on the configurations, except that no two tokens can share a vertex.

\begin{observation}\label{obs:simple-1-1-reconf}
Let $V_s$ and $V_t$ be two sets of vertices of size $k$ of a hypergraph $H$ of diameter $d$.
We can compute a reconfiguration sequence between $V_s$ and $V_t$ under \nkTJ{1}{1} that has at most $\mathcal{O}(n \cdot d)$ moves in $\mathcal{O}(n^2 m)$ time.
\end{observation}

For the proof, refer to \Cref{pf:simple-1-1-reconf}.

\begin{observation}\label{obs:reconf-keep-intersection}
Let $V_s$ and $V_t$ be two sets of vertices of size $k$ of a hypergraph $H$ of diameter $d$.
We can compute a reconfiguration sequence between $V_s$ and $V_t$ under \nkTJ{k}{1} that has at most $\mathcal{O}(n \cdot d)$ moves and has $V_s \cap V_t$ always occupied in $\mathcal{O}(n^2 m)$ time.
\end{observation}

For the proof, refer to \Cref{pf:reconf-keep-intersection}.

\subsection{Vertex cover and dominating set}
We present an efficient procedure that constructs a reconfiguration sequence between vertex covers in hypergraphs under \infTJ.
We consider the more general case of hypergraphs as that is useful in the later section on \probDS.
The idea is to find an intermediate vertex cover $V_m$, which can be reached by some tokens from both $V_s$ and $V_t$ in one move.
If this is possible, we can have some tokens occupy $V_m$ in all intermediate configurations and let the remaining tokens reconfigure while ensuring that $V_m$ is occupied, until we reach $V_t$.


\begin{lemma}\label{lem:move-to-common}
  Let $V_s, V_t$ be two vertex covers of a hypergraph $H = (V, F)$.
  Then we can compute in $\mathcal{O}(n^{3.5} + n^2 m)$ time a vertex cover $V_m$ such that both $B(V_s, V_m, H)$ and $B(V_t, V_m, H)$ have matchings saturating their copies of $V_m$.
\end{lemma}
\begin{proof}
  We say that a vertex cover $V_i$ is \emph{good} if both $B(V_s, V_i, H)$ and $B(V_t, V_i, H)$ have matchings that saturate the vertices corresponding to $V_i$.
  We present an algorithm which in each step either yields a good vertex cover $V_i$ or decreases the size of $V_i$.
  Thus, after at most $n$ steps, it yields a good vertex cover.

  Let $V_i$ be the vertex cover in $i$-th iteration, initially set $V_1 := V_s$.
  Now given $V_i$, check by \Cref{obs:computing-move} if both $B(V_s, V_i, H)$ and $B(V_t, V_i, H)$ have matchings that saturate their copies of $V_i$ and if they have, output $V_i$ as good.

  Otherwise, one of the bipartite graphs does not satisfy Hall's condition, without loss of generality, suppose it is $B' = B(V_s, V_i, H) = (V'_s \cup V'_i, E_B)$.
  Note that for any $A' \subseteq V'_i$, $A'$ corresponds to some $A \subseteq V_i$ and $N_{B'}(A')$ corresponds exactly to $N_H[A] \cap V_s$.
  Thus, by \Cref{thm:hall}, there is $A \subseteq V_i$ such that $|A| > |N_H[A] \cap V_s|$.
  Now we set $V_{i+1} \coloneqq (V_i \setminus A) \cup (N_H[A] \cap V_s)$, note that $|V_{i+1}| < |V_i|$.
  \begin{claim}
    $V_{i+1}$ is a vertex cover of $H$.
  \end{claim}
  \begin{proof}
    We claim that for each $e \in F$, if $e$ intersects $A$, then it intersects $N_H[A] \cap V_s$, from which the claim immediately follows.
    Suppose there is $e$ which intersects $A$ but does not intersect $N_H[A] \cap V_s$.
    Let $v \in e \cap A$ and let $v' \in e \cap V_s$, such $v'$ must exist as $V_s$ is a vertex cover.
    Then $dist(v, v') \leq 1$ and $v' \in N_H[A] \cap V_s \cap e$, a contradiction.
  \end{proof}
  In each iteration, the algorithm either yields a good $V_m$ or produces a vertex cover of smaller size, thus after at most $n$ iterations, it yields a good $V_m$.
  Finding a maximum matching in a bipartite graph and finding a Hall violator can be done in time $\mathcal{O}(n^{2.5} + n m)$, thus the total runtime is $\mathcal{O}(n^{3.5} + n^2 m)$.
\end{proof}

The following lemma shows that we can always find a $V'$ such that $V' \supseteq V_m$ and there is a move from $V_s$ to $V'$.

\begin{lemma}\label{lem:move-to-disjoint}
  Let $V_s$ and $V_m$ be sets of vertices of hypergraph $H$ such that there exists a matching in $B(V_s, V_m, H)$ saturating the copy of $V_m$.
  Then we can compute $V' \subseteq V(H)$ such that $V' \supseteq V_m$ along with a move from $V_s$ to $V'$ in $\mathcal{O}(n^2 m)$ time.
\end{lemma}
\begin{proof}
  First, let $M$ be the matching saturating the copy of $V_m$ in $B(V_s, V_m, H)$.
  Let $B'$ be the result of taking $B(V_s, V_m, H)$ and assigning to each edge a cost equal to the distance of the corresponding vertices in $G$.
  Let $M^*$ be the minimum cost matching saturating the copy of $V_m$ in $B'$, and let $V'_s \subseteq V_s$ be the vertices whose copy for $V_s$ is not in any edge of $M^*$.

  Let $V' = V'_s \cup V_m$, we need to show that $V'_s \cap V_m = \emptyset$.
  Suppose that there is $v \in V'_s \cap V_m$.
  The edge $e \in E(B')$ consisting of two copies of $v$ has weight~$0$ and $e \notin M^*$.
  As $M^*$ saturates the copy of $V_m$, there is $e' \in M^*$ of weight $1$ containing the copy of $v$ for $V_m$.
  Thus, $(M^* \setminus \{e'\}) \cup \{e\}$ is a matching saturating the copy of $V_m$ with a smaller weight, a contradiction.

  The move from $V_s$ to $V'$ has token on $V'_s$ stay and each token of $V_s \setminus V'_s$ will move along $M^*$ and thus reaching $V_m$.
  We can compute $B'$ in $\mathcal{O}(n^2 m)$ by finding the distances from each vertex of $H$.
  A minimum cost matching of $B'$ can be computed in time~$\mathcal{O}(n^3)$~\cite{west}.
\end{proof}

\begin{theorem}[\textbf{Vertex cover reachability}]\label{thm:reconf-hg}
  Let $H = (V, F)$ be a connected hypergraph of diameter $d$ and $V_s, V_t$ two vertex covers of $H$ with size $k$.
  We can compute a reconfiguration sequence under \infTJ that has at most $\mathcal{O}(n \cdot d)$ moves in $\mathcal{O}(n^{3.5} + n^2 m)$ time.
\end{theorem}
\begin{proof}
  First, compute $V_m$ by \Cref{lem:move-to-common} in time $\mathcal{O}(n^{3.5} + n^2 m)$ and then compute $V'_s$ and $V'_t$ with $V_m \subseteq V'_s, V_m \subseteq V'_t$ by \Cref{lem:move-to-disjoint}.
  We also have a move from $V_s$ to $V'_s$ and as each move is symmetric, also a move from $V'_t$ to $V_t$.
  It remains to obtain a reconfiguration sequence between $V'_s$ and $V'_t$.
  By \Cref{obs:reconf-keep-intersection}, we can compute a reconfiguration sequence from $V'_s$ to $V'_t$ while keeping $V'_s \cap V'_t \supseteq V'_m$ always occupied.
  Thus, every intermediate configuration contains $V_m$ or is $V_s$ or $V_t$ and therefore is a vertex cover.
\end{proof}

A similar idea is used in an algorithm that constructs a reconfiguration sequence between two dominating sets under \nknkTJ{k}{1}{k-2}{2}.

\begin{theorem}[\textbf{Dominating set reachability}]\label{prop:ds-connected}
  Let $G = (V, E)$ be a connected graph of diameter at most $d$ and $D_s$, $D_t$ two dominating sets of $G$ with size $k$.
  We can compute a reconfiguration sequence under \nknkTJ{k}{1}{k-2}{2} that has at most $\mathcal{O}(n \cdot d)$ moves in $\mathcal{O}(n^{3.5} + n^2 m)$ time.
\end{theorem}
\begin{proof}
  Let $H = (V, F)$ be a hypergraph, where $F$ is the set of all closed neighborhoods of $G$.
  Observe that $D$ is a dominating set of $G$ if and only if $D$ is a vertex cover of $H$.
  Furthermore, note that for $u, v \in V$ $dist_G(v, u) \leq 2$ holds if and only if there is $f \in F$ such that $u, v \in f$.


  Suppose that by \Cref{lem:move-to-common} applied on $H, D_s, D_t$ we computed $D_m$ and by \Cref{lem:move-to-disjoint} we computed $V'_s$ and $V'_t$.
  Let $B'_s = B(D_s, V'_s, H)$ and $B'_t = B(D_t, V'_t, H)$ and let $M_s, M_t$ be their matchings corresponding to the move from $D_s$ to $V'_s$ and the move from $D_t$ to $V'_t$, also obtained by \Cref{lem:move-to-disjoint}.
  We say that the distance of edge $\{u, v\}$ in $B'_s$ or $B'_t$ is the distance of the vertices corresponding to $u$ and $v$ in $G$.
  If both $M_s$ and $M_t$ have at most $k - 2$ edges with distance $2$, then they are valid moves under \nknkTJ{k}{1}{k-2}{2}.
  In the opposite case, we show how to modify each matching separately, so that at most $k - 2$ edges with distance $2$ are used.

  Without loss of generality, suppose that $M_s$ has at most $1$ edge with distance at most~$1$.
  Let $W_{\leq 1}$ be the edges of $B'_s$ with distance at most one.
  Observe that as $D_s, V'_s$ are dominating sets, for each $v \in D_s$, there is $u \in V'_s$ with $dist_G(v, u) \leq 1$.
  Hence, each vertex of $B'_s$ is incident to some edge in $W_{\leq 1}$.


  We show how to obtain a cycle of size at least $4$ consisting of alternating edges of $M_s$ and $W_{\leq 1}$, which is then used to exchange the edges of the matching with edges of $W_{\leq 1}$.
  If every vertex of $B'_s$ is incident to some edge in $W_{\leq 1} \setminus M_s$, we can easily find an alternating cycle by a standard argument. 
  If this is not the case, note that by assumption $|M_s \cap W_{\leq 1}| \leq 1$ and hence we reduce to the previous case by considering $B'_s$ with the vertices of $M_s \cap W_{\leq 1}$ deleted.


  Hence way obtain a cycle $C$ of size at least $4$ and with alternating edges of $M_s$ and $W_{\leq 1}$.
  We construct a new matching $M'_s \coloneqq M_s \triangle E(C)$ exchanging $|C| / 2 \geq 2$ edges for edges of distance $1$.
  The obtained $M'_s$ corresponds to a move from $D_s$ to $V'_s$ under \nknkTJ{k}{1}{k-2}{2}.
  The same argument can be applied for $D_t$ and $V'_t$.
  By \Cref{obs:reconf-keep-intersection}, we can obtain a reconfiguration sequence from $V'_s$ to $V'_t$ under \nkTJ{k}{1}, which has $V'_s \cap V'_t \supseteq D_m$ always occupied.
  Hence, the resulting sequence is valid under \nknkTJ{k}{1}{k-2}{2}.

\end{proof}

The following proposition shows that \nknkTJ{k}{1}{k-2}{2} is a minimal rule guaranteeing reachability.

\begin{figure}
  \centering
  \includegraphics{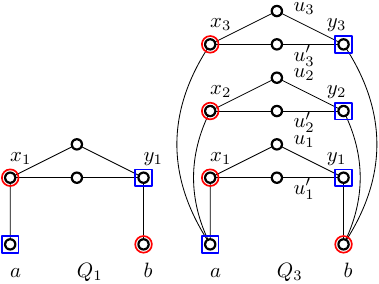}
  \caption{Construction of a graph, such that moving all tokens at once by distance at least $1$ and at least $k - 2$ of them by distance at least $2$
  is necessary to reconfigure between the blue and red dominating sets.}\label{fig:ds_not_connected}
\end{figure}

\begin{proposition}[\textbf{Dominating set unreachability}]\label{prop:ds-not-connected}
  For every $k \geq 3$, there exists a connected graph $G$ such that to reconfigure between two distinct minimum dominating sets of size $k$,
  we need to allow all tokens to move at once and at least $k - 2$ of them by distance at least $2$.
\end{proposition}
\begin{proof}
  Consider the graphs $Q_1$ and $Q_3$ shown in~\Cref{fig:ds_not_connected}.
  We can inductively construct $Q_i$ for $i \geq 2$ by taking a copy of $Q_{i-1}$ and $Q_1$ and identifying its vertices $a$ and $b$.

  We denote the size of a minimum dominating set of $Q_i$ as $\gamma(Q_i)$.
  First, we show that $\gamma(Q_i) \leq i + 1$ by constructing a dominating set on $i + 1$ vertices:
  observe that $\{a, y_1, y_2, \dots, y_i\}$ is a dominating set of size $i + 1$.
  On the other hand, from each set $S_j = \{x_j, u_j, u'_j, y_j\}$ we have to include at least one vertex.
  At the same time, no vertex from $S_j$ dominates all the other vertices in $S_j$ or any vertices from any $S_p$ for $p \neq j$, thus $\gamma(Q_i) \geq i + 1$.

  Now we claim that if $D$ is a minimum dominating set, then exactly one of $a, b$ is in $D$.
  Suppose that $a, b \in D$, then we need one additional vertex for each $u_i$, thus $|D| \geq 2 + i$, a contradiction.
  Suppose that $a, b \notin D$, then $D$ must contain at least $2$ vertices from each set $S_j$.
  Then we have $|D| \geq 2i$ and as $k \geq 3$ and hence $i \geq 2$, we have $|D| > i + 1$, a contradiction.

  Note that if $a \in D$ then $y_j \in D$ for every $j$.
  Thus, we have exactly two different minimum dominating sets, $D_a = \{a, y_1, y_2, \dots, y_i\}$ and $D_b = \{b, x_1, x_2, \dots, x_i\}$.
  Now observe that to move from $D_a$ to $D_b$, at least $k - 2$ tokens must move by distance at least $2$.
  Any matching between $D_a$ and $D_b$ such that matched tokens are at distance at most $2$ will match $k - 2$ tokens on some $x_j$ to some $y_p$.
  As $dist(x_j, y_p) \geq 2$ for all $j, q$, those tokens must move by distance at least $2$.
\end{proof}



\subsection{Independent sets}

The following will allow us to reduce reconfiguring independent sets to moving tokens to vertices of some shared maximal independent set.

\begin{lemma}\label{lem:indep-1-3}
  Let $G$ be a connected graph, $I^*$ a maximal independent set of $G$, and $I_s \subseteq I^*$, $I_t \subseteq I^*$ two independent sets of $G$ of size $k$.
  Then we can compute a reconfiguration sequence of length $\mathcal{O}(n^2)$ between $I_s$ and $I_t$ under \nkTJ{1}{3} in $\calO(n^2 m)$ time.
\end{lemma}
\begin{proof}
  Let $G' = (I^*, E' = \{ \{u, v\} \mid dist_G(u, v) \leq 3 \text{ and } u, v \in I^* \})$.
  Note that a reconfiguration sequence under $(1, 1)$-\TS on $G'$ corresponds to a sequence under $(1, 3)$-\TS on $G$ and any configuration (not necessarily an independent set) of tokens on $G'$ where no two tokens share a vertex is an independent set on $G$.

  First, observe that $G'$ is connected, to see that, suppose $G'$ has connected components $C_1, C_2, \dots, C_\ell$.
  Let $P = (p_1, \dots, p_q)$ be a shortest path in $G$ from $C_1$ to $C_2$ with $p_1 \in C_1$ and $p_q \in C_2$.
  Note that $q \geq 5$, $dist_G(p_3, C_1) = 2$, and $dist_G(p_3, C_i) \geq 2$ for all $i$, therefore $I^* \cup \{p_3\}$ is an independent set, a contradiction with $I^*$ being maximal.
  By \Cref{obs:simple-1-1-reconf}, we can reconfigure $I_s$ to $I_t$ under \nkTJ{1}{1} on $G'$ using $\mathcal{O}(n^2)$ moves.


  Note that at no point do two tokens share a vertex.
  Therefore, we construct a sequence that is a valid reconfiguration sequence on $G$ under $(1, 3)$-\TS.
\end{proof}

\begin{theorem}[\textbf{Independent set reachability}]\label{prop:indep-possible}
Let $G$ be a connected graph and $I_s$, $I_t$ two independent sets of size $k$ of $G$.
Then we can compute a reconfiguration sequence between $I_s$ and $I_t$ of length $\calO(n^2)$ under \nknkTJ{k}{1}{1}{3} in $\calO(n^{3.5} + n^2 m)$ time.
\end{theorem}
\begin{proof}
  We describe an algorithm that reconfigures $I_s$ to $I_t$.
  The idea is to first construct a maximal independent set $I^*$ such that $I_s$ and $I_t$ can move into a subset of $I^*$ in one move.
  The idea of construction of $I^*$ is similar to that of \Cref{lem:move-to-common}.

  We start by constructing an independent set $I_m$ such that $B(I_s, I_i, G)$ and $B(I_t, I_i, G)$ have a matching that saturates the copies of $I_s$ and $I_t$ respectively, in that case, we say that $I_m$ is \emph{good}.
  Initially set $I_1 := I_s$.
  Now in each step, we show that either $I_i$ is good or we can construct a larger independent set $I_{i+1}$.

  Suppose that $I_i$ is not good, then one of the bipartite graphs does not have a matching saturating the copies of $I_s$ or $I_t$, suppose it is $B' = B(I_s, I_i, G) = (I'_s \cup I'_i, E_B)$.
  By Hall's condition, there is $A' \subseteq I'_s$ such that $|A'| > |N_{B'}(A')|$, which in $G$ corresponds to some $A \subseteq I_s$ such that
  $|A| > |N_G[A] \cap I_i|$.

  We set $I_{i+1} := (I_{i} \setminus N_G[A]) \cup A$, observe that $|I_{i+1}| > |I_i|$ as $(I_i \setminus N_G[A]) \cap A = \emptyset$.
  We claim that $I_{i+1}$ is an independent set.
  To see that, suppose there is $\{u, v\} \in E$ such that $\{u, v\} \subseteq I_{i+1}$, note that $u, v \in A \subseteq I_s$ or $u, v \in I_i \setminus N_G[A] \subseteq I_i$ may not hold by assumption.
  In the remaining case $u \in I_i \setminus N_G[A], v \in A$ it holds $dist_G(u, v) \geq 2$.

  Thus, after at most $n$ iterations we obtain a good independent set $I_m$.
  Now we can greedily add vertices to $I_m$ until we obtain a maximal independent set $I^*$ with $I_m \subseteq I^*$.
  Given matchings $M_s, M_t$ saturating $I_s$ in $B(I_s, I_m, G)$ and $I_t$ in $B(I_t, I_m, G)$ respectively, we can move the tokens on $I_s$ along $M_s$ to some $I'_s \subseteq I_m$ in one move under $(k,1)$-\TS.
  Similarly, we can move from $I_t$ to $I'_t \subseteq I_m$ in one move.
  By \Cref{lem:indep-1-3}, we can now construct a reconfiguration sequence of length $\mathcal{O}(n^2)$ between $I'_s$ and $I'_t$ under \nkTJ{1}{3} in time $\calO(n m)$.
  Thus, we obtain a reconfiguration sequence between $I_s$ and $I_t$ under \nknkTJ{k}{1}{1}{3}.

  We check at most $n$ times whether $I_i$ is good, which by \Cref{obs:computing-move} takes $\calO(n^{2.5} + n m)$ time, thus the total runtime is $\calO(n^{3.5} + n^2 m)$.
\end{proof}

\begin{figure}
  \begin{center}
  \includegraphics{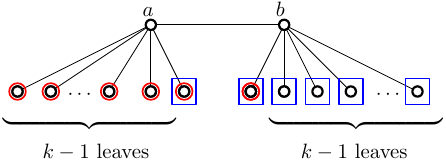}
  \caption{Sliding tokens to distance $3$ is necessary when reconfiguring between the red independent set and the blue independent set.}\label{fig:indep-not-connected}
  \end{center}
\end{figure}

We have shown in \Cref{prop:need-k-to-move} that at least \nkTJ{k}{1} is necessary to guarantee reachability.
Now, we show that \nkTJ{1}{3} is also necessary.

\begin{proposition}[\textbf{Independent set unreachability}]\label{prop:indep-not-connected}
  For every $k \geq 3$, there exists a connected graph $G$ such that to reconfigure between two distinct independent sets of size $k$,
  we need to allow at least one token to move by distance at least $3$.
\end{proposition}
\begin{proof}
  Consider the graph in \Cref{fig:indep-not-connected} with $k \geq 3$.
  Note that the tokens of the red independent set incident to vertex $a$ may not move, unless they jump to distance $3$.
  This holds for all the possible positions of the token incident to $b$, hence a move to distance $3$ is necessary.
\end{proof}

\section{Computational complexities}
In this section, we determine the complexities of deciding whether a reconfiguration sequence exists under \infTJ for the three studied problems, \probIS, \probVC, and \probDS.
It follows from \Cref{thm:reconf-hg} that the reachability problem for \probVC is in \Pclass for \infTJ, as the answer is always \YES.
In contrast, the problem of finding the \emph{shortest} reconfiguration sequence in this setting is shown to be \NP-complete.
We show that for \probIS and \probDS, the reachability question is \PSPACE-complete.

\begin{observation}\label{obs:in-pspace}
    \textsc{Reconfiguration of $\Pi$ under \infTJ}, where $k$ is the size of the reconfigured solution, is in \PSPACE for any problem $\Pi$ in \NP.
\end{observation}
\begin{proof}
    There is a simple \NPSPACE algorithm.
    For a given configuration $V_s$, we can nondeterministically enumerate all solutions of $\Pi$ and traverse into those that are reachable in a single move, which can be checked in polynomial time by~\Cref{obs:computing-move}.
    The computation stops after traversing $2^n$ configurations, as that is an upper bound for the length of the shortest reconfiguration sequence.
    By Savitch's Theorem~\cite{Savitch_1970}, this \NPSPACE algorithm can be converted into a \PSPACE one.
\end{proof}

\subsection{Independent sets}
In this section, we show that \probIS under \infTJ, where $k$ is the size of the reconfigured solution, is \PSPACE-complete.
This depends on the result of Belmonte et al., which shows that \probIS reconfiguration under \TS is \PSPACE-complete even on split graphs~\cite{Belmonte_2020}.
A graph is a \emph{split} graph if its vertices can be partitioned into a clique and an independent set.
We show that allowing multiple tokens to move at once provides no advantage on split graphs.

\begin{lemma}\label{lem:is-split}
  Let $G = (C \cup I, E)$ be a split graph with clique $C$ and independent set $I$ and $I_s$ and $I_t$ two independent sets of $G$ of size $k$.
  There is a reconfiguration sequence between $I_s$ and $I_t$ under \TS if and only if there is one under \infTJ.
\end{lemma}
\begin{proof}
  If there is a reconfiguration sequence under \TS, the same sequence is valid under \infTJ as well.

  Now suppose that there is a reconfiguration sequence between $I_s$ and $I_t$ under \infTJ.
  We show that such a sequence can be transformed so that each move consists of only one token slide and is thus valid under \TS.
  Note that in each configuration, $C$ can contain at most one token, and therefore we can move at most one token from $I$ to $C$ and move at most one token from $C$ to $I$ during one move.
  Therefore, each move consists of at most two token slides and each move with two token slides has one token move from $C$ and one move to $C$.

  Let $s = ((c_1, i_1), (i_2, c_2))$ be a move consisting of two token slides.
  Without loss of generality, let $c_1, c_2 \in C, i_1, i_2 \in I$.
  Note that if $i_1 = i_2$, we may reach the same configuration with a single token slide from $c_1$ to $c_2$.
  Hence, assume that $i_1 \neq i_2$.
  We claim that we can replace $s$ with two sequential moves $(c_1, i_1)$ and $(i_2, c_2)$ and the reconfiguration sequence remains valid.

  Let $I_i$ be the configuration in the reconfiguration sequence that is reached after applying $s$, that is $I_{i} = (I_{i-1} \setminus \{c_1, i_2\}) \cup \{i_1, c_2\}$.
  First note that $I' = (I_{i-1} \setminus \{c_1\}) \cup \{i_1\}$ (reached by sliding token from $c_1$ to $i_1$) is independent as only vertices of $I$ have tokens.
  Thus, in the reconfiguration sequence, we can place $I'$ between $I_{i-1}$ and $I_i$ to ensure that reaching $I_{i}$ from $I_{i-1}$ is done using single token slides.
  After modifying each move with two token slides, we arrive at a reconfiguration sequence under \TS.
\end{proof}

\begin{proposition}\label{thm:pts-hard}
  \probIS under \infTJ is \PSPACE-complete even on split graphs.
\end{proposition}
\begin{proof}
  Reconfiguration of \probIS under \TS is \PSPACE-complete~\cite{Belmonte_2020}.
  By \Cref{lem:is-split}, reconfiguration of \probIS under \infTJ is \PSPACE-complete as well.
\end{proof}

\subsection{Vertex covers}
In this section, we show that determining the length of the shortest reconfiguration sequence for \probVC under \infTJ is an \NP-complete problem.

\begin{theorem}\label{thm:vc-shortest-hard}
  Deciding if a reconfiguration sequence of length at most $\ell$ between two vertex covers exists under \infTJ is \NP-complete, even if $\ell \geq 2$ is fixed.
\end{theorem}

For the proof see \Cref{pf:vc-shortest-hard}

\subsection{Dominating sets}\label{section:dom}

For \probDS, we show that the problem remains \PSPACE-complete under \infTJ.
For that end, we define reconfiguration rule with additional allowed moves, which is useful in the subsequent hardness proof.
Its definition is in \Cref{subsec:appendix-dom}.

\begin{lemma}\label{lm:relaxed-seq}
  There exists a reconfiguration sequence between vertex covers $V_s$ and $V_t$ under relaxed token jumping if and only if there exists a reconfiguration sequence under \TJ.
\end{lemma}
For the proof, see \Cref{pf:relaxed-seq}.
This is used in the reduction from reconfiguration of \probVC under \TJ.

\begin{theorem}\label{thm:ds-reconf}
  Deciding if one dominating set can be reconfigured into another under \infTJ, where $k$ is the size of the given dominating sets, is \PSPACE-complete.
\end{theorem}
For the proof, see \Cref{pf:ds-reconf}

\section{Conclusion}
For \probDS and \probIS, we proved that \nknkTJ{k}{1}{k-2}{2} and \nknkTJ{k}{1}{1}{3} respectively ensure reachability, here $k$ denotes the size of the given reconfigured solutions.
If we decrease any of the parameters of these rules, we necessarily lose this guarantee.
A natural next step is to determine the complexity of finding a shortest reconfiguration sequence under these rules.

It is also an interesting question which graph classes allow weaker rules to guarantee reachability.
Furthermore, the investigated question is important for other graph problems as well, for instance, generalizations of \probDS and \probIS seem to be tractable with the tools we provide.

\subsubsection{\ackname}

\begin{credits}
J. M. Křišťan acknowledges the support of the Czech Science Foundation Grant No. \mbox{24-12046S}. This work was supported by the Grant Agency of the Czech Technical University in Prague, grant No. SGS23/205/OHK3/3T/18. J. Svoboda acknowledges the support of the ERC CoG 863818 (ForM-SMArt) grant.
\end{credits}

\bibliographystyle{splncs04}
\bibliography{main-parallel}

\begin{thebibliography}{10}
\providecommand{\url}[1]{\texttt{#1}}
\providecommand{\urlprefix}{URL }
\providecommand{\doi}[1]{https://doi.org/#1}

\bibitem{Bartier_2021}
Bartier, V., Bousquet, N., Dallard, C., Lomer, K., Mouawad, A.E.: On girth and
  the parameterized complexity of token sliding and token jumping. Algorithmica
   \textbf{83}(9),  2914--2951 (Jul 2021). \doi{10.1007/s00453-021-00848-1}

\bibitem{Belmonte_2020}
Belmonte, R., Kim, E.J., Lampis, M., Mitsou, V., Otachi, Y., Sikora, F.: Token
  sliding on split graphs. Theory of Computing Systems  \textbf{65}(4),
  662--686 (Mar 2020). \doi{10.1007/s00224-020-09967-8}

\bibitem{Bonamy_2013}
Bonamy, M., Bousquet, N.: Recoloring bounded treewidth graphs. Electronic Notes
  in Discrete Mathematics  \textbf{44},  257--262 (Nov 2013).
  \doi{10.1016/j.endm.2013.10.040}

\bibitem{Bonamy_2017}
Bonamy, M., Bousquet, N.: Token sliding on chordal graphs. In: Graph-Theoretic
  Concepts in Computer Science. pp. 127--139. Springer International Publishing
  (2017). \doi{10.1007/978-3-319-68705-6_10}

\bibitem{Bonamy_2021}
Bonamy, M., Dorbec, P., Ouvrard, P.: Dominating sets reconfiguration under
  token sliding. Discrete Applied Mathematics  \textbf{301},  6--18 (Oct 2021).
  \doi{10.1016/j.dam.2021.05.014}

\bibitem{Bonsma_2009}
Bonsma, P., Cereceda, L.: Finding paths between graph colourings:
  {PSPACE}-completeness and superpolynomial distances. Theoretical Computer
  Science  \textbf{410}(50),  5215--5226 (Nov 2009).
  \doi{10.1016/j.tcs.2009.08.023}

\bibitem{Bonsma_2014}
Bonsma, P., Kamiński, M., Wrochna, M.: Reconfiguring independent sets in
  claw-free graphs. In: Algorithm Theory – SWAT 2014. pp. 86--97. Springer
  International Publishing (2014). \doi{10.1007/978-3-319-08404-6_8}

\bibitem{Bonsma_2014_2}
Bonsma, P., Mouawad, A.E., Nishimura, N., Raman, V.: The complexity of bounded
  length graph recoloring and {CSP} reconfiguration. In: Parameterized and
  Exact Computation. pp. 110--121. Springer International Publishing (2014).
  \doi{10.1007/978-3-319-13524-3_10}

\bibitem{Bousquet_2021}
Bousquet, N., Joffard, A.: {TS}-reconfiguration of dominating sets in circle
  and circular-arc graphs. In: Fundamentals of Computation Theory. pp.
  114--134. Springer International Publishing (2021).
  \doi{10.1007/978-3-030-86593-1_8}

\bibitem{Censor_Hillel_2020}
Censor-Hillel, K., Rabie, M.: Distributed reconfiguration of maximal
  independent sets. Journal of Computer and System Sciences  \textbf{112},
  85--96 (Sep 2020). \doi{10.1016/j.jcss.2020.03.003}

\bibitem{Cereceda_2010}
Cereceda, L., van~den Heuvel, J., Johnson, M.: Finding paths between
  3-colorings. Journal of Graph Theory  \textbf{67}(1),  69--82 (Dec 2010).
  \doi{10.1002/jgt.20514}

\bibitem{Charrier_2020}
Charrier, T., Queffelec, A., Sankur, O., Schwarzentruber, F.: Complexity of
  planning for connected agents. Autonomous Agents and Multi-Agent Systems
  \textbf{34}(2) (Jun 2020). \doi{10.1007/s10458-020-09468-5}

\bibitem{Demaine_2015}
Demaine, E.D., Demaine, M.L., Fox-Epstein, E., Hoang, D.A., Ito, T., Ono, H.,
  Otachi, Y., Uehara, R., Yamada, T.: Linear-time algorithm for sliding tokens
  on trees. Theoretical Computer Science  \textbf{600},  132--142 (Oct 2015).
  \doi{10.1016/j.tcs.2015.07.037}

\bibitem{Gajjar_2022}
Gajjar, K., Jha, A.V., Kumar, M., Lahiri, A.: Reconfiguring shortest paths in
  graphs. Proceedings of the AAAI Conference on Artificial Intelligence
  \textbf{36}(9),  9758--9766 (Jun 2022). \doi{10.1609/aaai.v36i9.21211}

\bibitem{Gan_2019}
Gan, J., Suksompong, W., Voudouris, A.A.: Envy-freeness in house allocation
  problems. Mathematical Social Sciences  \textbf{101},  104--106 (Sep 2019).
  \doi{10.1016/j.mathsocsci.2019.07.005}

\bibitem{Gopalan_2006}
Gopalan, P., Kolaitis, P.G., Maneva, E.N., Papadimitriou, C.H.: The
  connectivity of boolean satisfiability: Computational and structural
  dichotomies. In: Automata, Languages and Programming. pp. 346--357. Springer
  Berlin Heidelberg (2006). \doi{10.1007/11786986_31}

\bibitem{Haddadan_2016}
Haddadan, A., Ito, T., Mouawad, A.E., Nishimura, N., Ono, H., Suzuki, A.,
  Tebbal, Y.: The complexity of dominating set reconfiguration. Theoretical
  Computer Science  \textbf{651},  37--49 (Oct 2016).
  \doi{10.1016/j.tcs.2016.08.016}

\bibitem{Hall_1935}
Hall, P.: On representatives of subsets. Journal of the London Mathematical
  Society  \textbf{s1-10}(1),  26--30 (Jan 1935).
  \doi{10.1112/jlms/s1-10.37.26}

\bibitem{Hoang_2019}
Hoang, D.A., Khorramian, A., Uehara, R.: Shortest reconfiguration sequence for
  sliding tokens on spiders. In: Algorithms and Complexity. pp. 262--273.
  Springer International Publishing (2019). \doi{10.1007/978-3-030-17402-6_22}

\bibitem{ITO_2016}
Ito, T., Nooka, H., Zhou, X.: Reconfiguration of vertex covers in a graph.
  IEICE Transactions on Information and Systems  \textbf{E99.D}(3),  598--606
  (2016). \doi{10.1587/transinf.2015fcp0010}

\bibitem{Kami_ski_2011}
Kamiński, M., Medvedev, P., Milanič, M.: Shortest paths between shortest
  paths. Theoretical Computer Science  \textbf{412}(39),  5205--5210 (Sep
  2011). \doi{10.1016/j.tcs.2011.05.021}

\bibitem{Kami_ski_2012}
Kamiński, M., Medvedev, P., Milanič, M.: Complexity of independent set
  reconfigurability problems. Theoretical Computer Science  \textbf{439},
  9--15 (Jun 2012). \doi{10.1016/j.tcs.2012.03.004}

\bibitem{Klostermeyer_2016}
Klostermeyer, W., Mynhardt, C.: Protecting a graph with mobile guards.
  Applicable Analysis and Discrete Mathematics  \textbf{10}(1),  1--29 (2016).
  \doi{10.2298/aadm151109021k}

\bibitem{K_i_an_2023}
Křišťan, J.M., Svoboda, J.: Shortest dominating set reconfiguration under
  token sliding. In: Fundamentals of Computation Theory. pp. 333--347. Springer
  Nature Switzerland (2023). \doi{10.1007/978-3-031-43587-4_24}

\bibitem{Lokshtanov_2019}
Lokshtanov, D., Mouawad, A.E.: The complexity of independent set
  reconfiguration on bipartite graphs. ACM Transactions on Algorithms
  \textbf{15}(1),  1--19 (Oct 2018). \doi{10.1145/3280825}

\bibitem{Lokshtanov_2018}
Lokshtanov, D., Mouawad, A.E., Panolan, F., Ramanujan, M., Saurabh, S.:
  Reconfiguration on sparse graphs. Journal of Computer and System Sciences
  \textbf{95},  122--131 (Aug 2018). \doi{10.1016/j.jcss.2018.02.004}

\bibitem{Mouawad_2018}
Mouawad, A., Nishimura, N., Raman, V., Siebertz, S.: Vertex cover
  reconfiguration and beyond. Algorithms  \textbf{11}(2), ~20 (Feb 2018).
  \doi{10.3390/a11020020}

\bibitem{Mouawad_2017}
Mouawad, A.E., Nishimura, N., Pathak, V., Raman, V.: Shortest reconfiguration
  paths in the solution space of boolean formulas. SIAM Journal on Discrete
  Mathematics  \textbf{31}(3),  2185--2200 (Jan 2017). \doi{10.1137/16m1065288}

\bibitem{Nishimura_2018}
Nishimura, N.: Introduction to reconfiguration. Algorithms  \textbf{11}(4), ~52
  (Apr 2018). \doi{10.3390/a11040052}

\bibitem{Savitch_1970}
Savitch, W.J.: Relationships between nondeterministic and deterministic tape
  complexities. Journal of Computer and System Sciences  \textbf{4}(2),
  177--192 (Apr 1970). \doi{10.1016/s0022-0000(70)80006-x}

\bibitem{Stern_2019}
Stern, R.: Multi-agent path finding – an overview. In: Artificial
  Intelligence. pp. 96--115. Springer International Publishing (2019).
  \doi{10.1007/978-3-030-33274-7_6}

\bibitem{Tateo_2018}
Tateo, D., Banfi, J., Riva, A., Amigoni, F., Bonarini, A.: Multiagent connected
  path planning: {PSPACE}-completeness and how to deal with it. Proceedings of
  the AAAI Conference on Artificial Intelligence  \textbf{32}(1) (Apr 2018).
  \doi{10.1609/aaai.v32i1.11587}

\bibitem{west}
West, D.B.: Introduction to graph theory. Prentice Hall, Upper Saddle River, NJ
  (2001)

\bibitem{Wrochna_2018}
Wrochna, M.: Reconfiguration in bounded bandwidth and tree-depth. Journal of
  Computer and System Sciences  \textbf{93},  1--10 (May 2018).
  \doi{10.1016/j.jcss.2017.11.003}

\bibitem{Yamada_2021}
Yamada, T., Uehara, R.: Shortest reconfiguration of sliding tokens on
  subclasses of interval graphs. Theoretical Computer Science  \textbf{863},
  53--68 (Apr 2021). \doi{10.1016/j.tcs.2021.02.019}

\end{thebibliography}

\appendix
\section{Omitted proofs}
\subsection{Proof of \Cref{obs:simple-1-1-reconf}}
\label{pf:simple-1-1-reconf}
\begin{proof}
    We can reconfigure $V_s$ to $V_t$ under \nkTJ{1}{1} on $H$ as follows.
    In each iteration, we increase the size of $V_s \cap V_t$.

    First, we pick $s \in V_s \setminus V_t, t \in V_t \setminus V_s$ and over several token slides, ensure that the token on $s$ is moved to $t$.
    Let $P = p_1 = s, p_2, \dots, p_\ell = t$ be the vertices of a shortest path between $s$ and $t$.
    Suppose the token is currently placed on $p_i$, if $p_{i+1}$ has no token, then we may directly move the token from $p_i$ to $p_{i+1}$.
    Otherwise, suppose $p_{i+1}, p_{i+2}, \dots, p_{i+j}$ are occupied by tokens and $p_{i+j+1}$ is not occupied, such $p_{i+j+1}$ must exist as $t$ is not occupied by assumption.
    We may reach the same configuration corresponding to jumping with the token from $p_{i}$ to $p_{i+j+1}$ by performing a sequence of token slides from $p_{i+j-q}$ to $p_{i+j+1-q}$, for all $0 \leq q \leq i$.
    Note that each such move is valid as for $q = 0$, we have that $p_{i+j+1}$ is unoccupied by assumption and for $q \geq 1$, we have that $p_{i+j+1-q}$ was left unoccupied by the preceding move.

    A simple implementation of this process will in each iteration pick two arbitrary vertices in $V_s \setminus V_t$ and $V_t \setminus V_s$ and find a shortest path between them, which can be done in $\calO(n m)$ on hypergraphs.
    This is done at most $n$ times, thus the total runtime is $\calO(n^2 + m)$.
    The total number of moves in each iteration is the length of the shortest path, hence we perform at most $\mathcal{O}(n \cdot d)$ moves, where $d$ is the diameter of $H$.
\end{proof}

\subsection{Proof of \Cref{obs:reconf-keep-intersection}}
\label{pf:reconf-keep-intersection}
\begin{proof}
    Let $H'$ be a hypergraph so that $V(H') = V(H) \setminus (V_s \cap V_t)$ and $E(H')$ is constructed by taking $E(H)$ and adding to $H'$ the set of edges $E'$ connecting by an edge every pair of vertices reachable by a path with internal vertices in $V_s \cap V_t$.
    By \Cref{obs:simple-1-1-reconf}, we can obtain a reconfiguration sequence under \nkTJ{1}{1} between $V_s \setminus V_t$ and $V_t \setminus V_s$ on $H'$.
    This can be transformed to a reconfiguration sequence under \nkTJ{k}{1} between $V_s$ and $V_t$ on $H$ by replacing any move of a token over an edge of $E' = \{u, v\}$ by a move consisting of sliding tokens along the path connecting $u$ with $v$ through $V_s \cap V_t$.
    Note that this ensures that $V_s \cap V_t$ is always occupied.
\end{proof}

\subsection{Proof of \Cref{thm:vc-shortest-hard}}
\label{pf:vc-shortest-hard}
\begin{proof}
    \begin{figure}
        \centering
        \includegraphics[page=2]{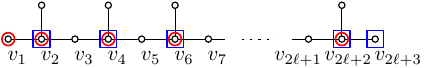}
        \caption{Gadget $T^\ell$. The vertices with red circles are in $V^\ell_s$ and the vertices with blue squares are in $V^\ell_t$.}\label{fig:vc-t}
    \end{figure}

    First we show that the problem is in \NP.
    Let $G$ be the input graph and $d$ its diameter.
    For $\ell > C \cdot n \cdot d$, the answer is always \YES for some fixed $C$ by \Cref{thm:reconf-hg}.
    Also, any reconfiguration sequence of length at most $\ell = \mathcal{O}(n \cdot d)$ can be verified in polynomial time using \Cref{obs:computing-move}.

    We follow with a reduction from \probVC.
    First let us define graph $T^\ell$, see \Cref{fig:vc-t} for an illustration.
    To construct $T^\ell$, first we create a path of length $2\ell + 3$ on vertices $v_1, v_2, \cdots, v_{2 \ell + 3}$.
    We connect a new vertex $u_{2k}$ with $v_{2k}$ for all integer $k$ such that $2 \leq 2k \leq 2 \ell + 2$.
    On $T^\ell$ we will later use as the starting configuration $V^{\ell}_s = \{v_1, v_2, v_4, v_6, \dots, v_{2\ell+2}\}$ and as the ending configuration $V^{\ell}_t = \{v_2, v_4, v_6, \dots, v_{2\ell+2}, v_{2\ell+3}\}$.

    \begin{claim}\label{clm:tl-reconf}
    Reconfiguring $V^{\ell}_s$ to $V^{\ell}_t$ on $T^\ell$ under \infTJ requires at least $\ell + 1$ turns.
    \end{claim}
    \begin{proof}
        The proof is by induction on $\ell$.
        First note that for reconfiguring $V^1_s$ to $V^1_t$ on $T^1$ we need at least $2$ moves.

        Now suppose that $\ell \geq 2$.
        Observe that exactly two different configurations are reachable from $V^\ell_s$ in one move:
        $S_1 = \{u_2, v_2, v_4, \dots, v_{2\ell + 2}\}$
        and $S_2 = \{v_2, v_3, v_4, \dots, v_{2\ell + 2}\}$.
        Only from $S_2$ can we reach any other configuration than $V^{\ell}_s, S_1, S_2$, thus that is the only possible next configuration in a shortest reconfiguration sequence.

        Note that $T \setminus \{v_1, v_2, u_2\}$ is isomorphic to $T^{\ell - 1}$ and $S_2 \setminus \{v_2\}$ corresponds to $V^{\ell - 1}_s$ on $T^{\ell - 1}$.
        The token on $v_2$ can not move in subsequent moves, otherwise, the edge connecting to its leaf would be uncovered.
        Thus by the induction hypothesis, we need at least $\ell$ additional moves to reconfigure from $S_2$ to $V^\ell_t$.
    \end{proof}

    The idea of the subsequent reduction is that if we connect $v_4$ to a new vertex, which can supply an additional token to $v_4$, we can make the reconfiguration sequence shorter by one move.
    This additional token will be available to every copy of $T^\ell$ only if a given graph a has small enough vertex cover.

    Now, we construct the instance $(G', V_s, V_t)$ of \textsc{Shortest Reconfiguration of} \probVC \textsc{under} \infTJ given an instance of \probVC $(G, k)$ and an integer $\ell \geq 2$.
    To construct $G'$, we first add $n - k$ copies of the gadget $T^\ell$, say $T_1, \dots, T_{k}$.
    We will distinguish the vertices of $T_i$ by superscript, meaning that path of $T_i$ consists of $v^i_1, \dots, v^i_{2\ell + 3}$.
    For each $T_i$, we connect $v^i_4$ with all vertices in $G$.
    This finishes the construction of $G'$.
    We set
    \[
        V_s \coloneqq \{v^i_{k} \mid i \in \{1, \dots, n - k\}, k \in \{1, 2, 4, 6, \dots, 2\ell + 2\}\} \cup V(G)
    \]
    \[
        V_t \coloneqq (V_s \setminus \{v^1_1, v^2_1, \dots, v^{n-k}_1\}) \cup \{v^1_{2\ell+3}, v^2_{2\ell+3}, \dots, v^{n-k}_{2\ell+3}\}
    \]
    \begin{claim}
        If $G$ has a vertex cover of size $k$ then there is a reconfiguration sequence from $V_s$ to $V_t$ using at most $\ell$ moves.
    \end{claim}
    \begin{proof}
        Let $V_c$ be a vertex cover of $G$ of size $k$ and $\overline{V_c} = V(G) \setminus V_c = \{c_1, \dots, c_{n-k}\}$.
        We construct the reconfiguration sequence as follows.
        The first move is\\
        $(v^i_1, v^i_2), (v^i_2, v^i_3), (c_i, v^i_4), (v^i_4, v^i_5)$ for all $i \in \{1, \dots, n - k\}$ simultaneously.
        The second move is $(v^i_3, v^i_4), (v^i_4, c_i), (v^i_5, v^i_6), (v^i_6, v^i_7)$ for all $i \in \{1, \dots, n - k\}$ simultaneously.
        Now for $j \in \{3, \dots, \ell\}$, the $j$-th move is $(v^i_{2j+1}, v^i_{2j+2}), (v^i_{2j+2}, v^i_{2j+3})$ for all $i \in \{1, \dots, n - k\}$.
        Note that each resulting state is a vertex cover, as $V_c$ always remains on $G$ and for all $i, j$ we have a token on $v^i_{j}$ for all $j \in \{2, 4, \dots, 2\ell + 2\}$.
        Thus, we reach $V_t$ in exactly $\ell$ moves.
    \end{proof}

    \begin{claim}
        If $G$ has no vertex cover on $k$ vertices, then any reconfiguration sequence between $V_s$ and $V_t$ has at least $\ell + 1$ moves.
    \end{claim}
    \begin{proof}
        As at least $k+1$ tokens must always be on $G$, there is a $T_q$ such that it has the same number of tokens before and after the first move.
        As observed in the proof of the previous claim, the only viable move on $T_q$ is $(v^q_1, v^q_2), (v^q_2, v^q_3)$.
        Then the token slides between $v^q_4$ and $G$ can not decrease the subsequent number of moves required to reconfigure the tokens on $T_q$ and by the previous claim we require $\ell + 1$ moves in total.
    \end{proof}
    Hence, $G$ has a vertex cover of size $k$ if and only if there is a reconfiguration sequence between $V_s$ and $V_t$ of length at most $\ell$, which completes the proof.
\end{proof}

\subsection{Proofs of the \Cref{section:dom}}
\label{subsec:appendix-dom}

We first consider a reconfiguration rule with additional allowed moves, which is useful in the subsequent hardness proof.

\begin{definition}
    Reconfiguration under \emph{relaxed token jumping} allows moves to remove a token, add a token, or jump with a token from one vertex to another.
    For each intermediate configuration $V_i$, it must hold $|V_i| \leq |V_s|$, where $V_s$ is the starting configuration.
    In one turn, we may perform at most $|V_s| - |V_i|$ additions and at most $|V_s| - |V_i| + 1$ jumps and additions in total.
    Furthermore, it is allowed to jump or add a token to an already occupied vertex.
\end{definition}

Formally, the relaxed sequence is a sequence of \emph{moves}, each consisting of a sequence of \emph{steps}.
Thus, we denote a relaxed token jumping sequence as \[((o_{1,1}, \dots, o_{1,\ell_1}), (o_{2,1}, \dots, o_{2,\ell_2}), \dots, (o_{p,1}, \dots, o_{p,\ell_p}))\] where $p$ is the total number of moves, $\ell_i$ is the number of steps in the $i$-th move, and $o_{i,j}$ is one of the steps done in parallel in the $i$-th move.
Each $o_{i,j} \in \{R(v), A(v), J(u, v)\}$, where $R(v)$ denotes a removal of a token on vertex $v$, $A(v)$ denotes an addition of a token to vertex $v$ and $J(u, v)$ denotes a jump of a token from vertex $u$ to vertex $v$.

\subsubsection{Proof of \Cref{lm:relaxed-seq}}
\label{pf:relaxed-seq}
\begin{proof}
    Observe that a reconfiguration sequence under \TJ is valid under relaxed token jumping.
    The rest of the proof shows how to transform a relaxed token jumping sequence into a \TJ sequence.
    We will show a set of operations that modify a valid relaxed sequence into another valid relaxed sequence, eventually reaching a form of only singleton jumps.
    Suppose we process the moves of a reconfiguration sequence $S$ under relaxed token jumping between $V_s$ into $V_t$ one by one.
    Assume that $S$ contains moves with more than one step.
    First, note that if a move contains a pair of an addition $A(v)$ and a removal $R(u)$, we may replace it with a jump $J(v, u)$.
    Hence, we may assume that all moves do not contain additions and removals at the same time.

    Now, note that whenever a move $m_i$ contains more than one step and one of them is a removal, say $R(v)$, we can remove $R(v)$ from $m_i$ and place a new move $m' = (R(v))$ directly after $m_i$.
    Similarly, if a move $m_i$ contains an addition, say $A(v)$, we may remove $A(v)$ from $m_i$ and place a new move $m' = (A(v))$ directly before $m_i$.
    The resulting configurations remain supersets of what they were before and the following moves remain valid.

    Thus, now we assume that removals and additions appear only in moves with one step.
    Now, any move with multiple steps consists of only jumps, note that all such jumps $J(u, u')$ can be replaced by a pair $R(u), A(u')$ and the previous reduction can be applied again.
    Hence, now all moves contain only one step.



    Now, if any removals remain, let $m_i = (R(v))$ be such move with maximum $i$.
    Let $m_j$ be the move with $j > i$ with minimum $j$ which would not be allowed without $m_i$ and hence contains an addition.
    First, note that we may place $m_i$ so that it directly precedes $m_j$ and the reconfiguration sequence remains valid by the choice of $j$ and the fact that each
    configuration remains a superset of what it was before the change.
    As $m_j$ contains an addition $A(u)$, we may replace $m_{j-1} = (R(v))$ with $m'_{j-1} = (J(v, u))$.
    Note that the reconfiguration sequence remains valid for all subsequence moves.
    After repeating this procedure, the resulting sequence consists of only single jumps, but we may have jumped with a token to an occupied vertex.

    Let $m_i = (J(u, v))$ such that in $V_i$, $v$ already has a token and such that $i$ is minimum possible.
    There must be $m_j = (J(v, y))$ such that $i < j$, as $V_t$ has no two tokens sharing a vertex.
    Assume $m_j$ is such with minimum $j$.
    Let $x$ be some vertex without a token in $V_i$.
    Then replace $m_i$ with $m'_i = (J(u, x))$ and $m_j$ with $m_j = (J(x, y))$.
    Note that each such replacement increases the index of the first move which places a token on an already occupied vertex,
    thus eventually all jumps are only to unoccupied vertices.
    Thus the resulting reconfiguration sequence is valid under \TJ.
\end{proof}

\subsubsection{Proof of \Cref{thm:ds-reconf}}
\label{pf:ds-reconf}
\begin{proof}
    The problem is in \PSPACE by \Cref{obs:in-pspace}.

    We show PSPACE-hardness by a reduction from reconfiguring \probIS under \TJ, which is known to be PSPACE-complete~\cite{Kami_ski_2012}.
    Note that it is equivalent to reconfiguring \probVC under \TJ (shortened to VC-TJ), as the complement of an independent set is always a vertex cover.
    Furthermore, the jump from $u$ to $v$ in an independent set is equivalent to a jump from $v$ to $u$ in the corresponding vertex cover.
    Thus, the two problems are equivalent and we continue our reduction from the reconfiguring \probVC under \TJ.
    The construction used in the reduction is similar to the standard reduction of \probVC to \probDS in bipartite graphs.

    Let $(G, V_s, V_t)$ be the VC-TJ instance on the input with $k = |V_s| = |V_t|$.
    We construct the instance of reconfiguration of \probDS under \infTJ (shortened to DS-$k$TJ) $(G', D_s, D_t)$ as follows.
    Let $G' \coloneqq (V(G) \cup E_1 \cup \dots \cup E_{k+1} \cup \{x, y\}, E')$, where each $E_i$ for $1 \leq i \leq k + 1$ is a distinct copy of $E(G)$.
    For each $e_i \in E_i$ corresponding to $\{u, v\} \in E(G)$, connect $u$ and $v$ with $e_i$ in $G'$.
    Also, connect $x$ with all vertices in $V(G) \cup \{y\}$.
    Now set $D_s \coloneqq V_s \cup \{x\}$ and $D_t \coloneqq V_t \cup \{x\}$.

    We will show that there is a reconfiguration sequence for $(G, V_s, V_t)$ under VC-TJ if and only if there is a reconfiguration sequence for $(G', D_s, D_t)$ under DS-$k$TJ.
    Suppose there is a sequence of states $V_s = V_1, V_2, \cdots, V_t = V_k$ reconfiguring under VC-TJ.
    We show how to construct the sequence for $(G', D_s, D_t)$.
    Each jump from $u$ to $v$ is translated to two parallel moves, one from $u$ to $x$ and another from $x$ to $v$.
    Now each resulting configuration of tokens on $G'$ is a dominating set, as each $v \in V(G) \cup \{y\}$ is dominated from $x$ and each $e \in E(G)$ is incident to some $v \in V(G)$ with a token, as each corresponding configuration on $G$ is a vertex cover.
    Thus, we can reconfigure $D_s$ into $D_t$.

    Now suppose there is a sequence $S'$ of configurations $D_s = D_1, D_2, \cdots, D_t = D_k$ reconfiguring under DS-$k$TJ.
    Note we can assume that $x \in D_i$ for all $i$ as a token must be on $x$ and $y$ to dominate $y$ and $N_{G'}[y] \subseteq N_{G'}[x]$.
    Furthermore, note that each $D_i$ must have a token on one vertex of each $e = \{u, v\} \in E(G)$, otherwise a copy of $e_i$ in some $E_i$ would be left undominated.
    Any move of a token to such $e_i \in E_i$ in $S'$ can be rerouted to $x$ as its neighborhood is a strict superset, $u$ and $v$ is already dominated from $x$, and one of $u$ and $v$ must have a token and hence $e_i$ is already dominated.
    Thus, we will assume that moves of tokens from $V(G)$ are only to $x$.
    We show how to construct the sequence $S$ for $(G, V_s, V_t)$ under relaxed VC-TJ.

    Each step in $S$ is constructed from one step in $S'$.
    Whenever in $S'$ a token moves from $v$ to $x$ and in parallel another token moves from $x$ to $u$, we translate that into a move in $S$ with a jump from $v$ to $u$.
    We pick such a matching of parallel moves in each step of $S'$ from $x$ and to $x$ arbitrarily.
    If, after considering the matching, a step in $S'$ has any moves from $x$ to $v_1, \cdots, v_p$ left, we translate that to additions of a token and placing it on $v_1, \cdots, v_p$ in $S$.
    Similarly, if after the matching, the current step has any moves to $x$ from $v_1, \cdots, v_p$ left, we translate that to removal of tokens on $v_1, \cdots, v_p$ in $S$.
    Note that each state reached using $S$ is still a vertex cover of $G$: the corresponding dominating set reached using $S'$ must dominate the vertices of all $E_i$ of $G'$ from $V(G)$.
    Furthermore, to perform $p$ jumps and additions in parallel, there must be at least $p$ tokens on $x$, thus on $V(G)$ is at most $k - p + 1$.
    Also, there will always be at most $k$ tokens on $V(G)$.
    Thus, $S$ is a valid reconfiguration between $V_s$ and $V_t$ under relaxed token jumping and by \Cref{lm:relaxed-seq}, there exists a reconfiguration sequence under \TJ for $(G, V_s, V_t)$.
    This concludes the proof.
\end{proof}

\end{document}